# APPROXIMATION OF PERIODIC FUNCTION BY FEJER SUM AND DE LA VALLEE POUSSIN SUM

## Mikhael Shahoud[*]

**Abstract**: in this paper there are several results, we prove approximation of periodic function by Fejer means and De La Vallee Poussin means in Lebesgue spaces $L_{2\pi}^{p}$ the estimates are given in terms of function for $L_2$ and in terms of second continuity modulus.

1. **Introduction:** for Hilbert space $L_2$ consists of all $2\pi$ - periodic $2^{th}$ power Lebesgue integrable functions for R with the norm

$$\|f\|_2 = \left\{ \frac{1}{2\pi} \int_{-\pi}^{\pi} |f(x)|^2 dx \right\}^{\frac{1}{2}}$$

Consider its trigonometric Fourier series for function $f(x)$

$$\frac{a_0}{2} + \sum_{k=1}^{\infty} \alpha_k \cos kx + \beta_k \sin kx \qquad (1)$$

And $S_n$ be the $n^{th}$ partial sum of Fourier series

$$S_n = \frac{a_0}{2} + \sum_{k=1}^{N} \rho_k \cos(kx + \varphi_k) \qquad (2)$$

And Fejer sum is one of the arithmetic means of the partial sums of a Fourier series for the first $n^{th}$ terms.

$$\sigma_{n-1}(f) = \frac{S_0 + S_1 + ... + S_{n-1}}{n}$$

**Lemma (1):** Fejer sums can be written in the form of a trigonometric polynomial of degree no greater than (n-1). (see [2]).
**Proof:**

We have from equation (2):

*\* Assistant professor of Mathematics, Wadi international University*

$$S_0 = \frac{a_0}{2}$$

$$S_1 = \frac{a_0}{2} + \rho_1 \cos(x + \varphi_1)$$

$$S_2 = \frac{a_0}{2} + \rho_1 \cos(x + \varphi_1) + \rho_2 \cos(2x + \varphi_2)$$

$$S_n = \frac{a_0}{2} + \rho_1 \cos(x + \varphi_1) + \rho_2 \cos(2x + \varphi_2) + \ldots + \rho_{n-1} \cos((n-1)x + \varphi_{n-1})$$

By sum and dived by n:

$$\sigma_{n-1}(f) = \frac{S_0 + S_1 + \ldots + S_{n-1}}{n}$$

$$\Rightarrow \sigma_{n-1}(f) = \frac{a_0}{2} + \sum_{k=1}^{n-1}\left(1 - \frac{k}{n}\right)\rho_k \cos(kx + \varphi_k)$$

**Remark (1):**

For every $f(x) \in L_2$, the series $\sum_{k=1}^{\infty} \rho_k^2$ will be convergent (see[4])

Moreover the terms of this series are non-negative real numbers, non increase from term $k = n_0$

Therefore

$$\rho_p^2 \leq \rho_q^2 \quad \forall p \quad , \quad (p > q \geq n_0)$$

We obtain

$$\frac{\rho_p^2}{p} \leq \frac{\rho_q^2}{q} \quad \forall p \quad , \quad (p > q \geq n_0)$$

**Remark (2):**

From Fourier series for function (1): $f(x) \approx \frac{a_0}{2} + \sum_{k=1}^{\infty} \alpha_k \cos kx + \beta_k \sin kx$

Therefore

$$f(x) - \sigma_{n-1}(f) \approx \sum_{k=1}^{n-1} \frac{k}{n} \rho_k \cos(kx + \varphi_k) + \sum_{k=n}^{\infty} \rho_k \cos(kx + \varphi_k)$$

From Parseval equality

$$\|f(x)-\sigma_{n-1}(f)\|^2 \approx \sum_{k=1}^{n-1}\left(\frac{k}{n}\right)^2 \rho_k^2 + \sum_{k=n}^{\infty}\rho_k^2 \qquad (3).$$

And De La Vallee Poussin sum is one of the arithmetic means of the partial sums of a Fourier series from $m^{th}$ term to $(n-1)^{th}$ term.

$$V_m^{n-1}(f) = \frac{S_m + S_{m+1} + \ldots\ldots + S_{n-1}}{n-m}$$

**Lemma (2):** De La Vallee Poussin sums can be written in the form of a trigonometric polynomial of degree no greater than (n-1). (see[3])
**Proof:**

We have from equation (2):

$$S_m = \frac{a_0}{2} + \sum_{k=1}^{m}\rho_k \cos(kx+\varphi_k)$$
$$S_{m+1} = S_m + \rho_{m+1}\cos[(m+1)x+\varphi_{m+1}]$$
$$S_{m+2} = S_m + \rho_{m+1}\cos[(m+1)x+\varphi_{m+1}] + \rho_{m+2}\cos[(m+2)x+\varphi_{m+2}]$$
.
.
$$S_{n-1} = S_m + \rho_{m+1}\cos[(m+1)x+\varphi_{m+1}] + \rho_{m+2}\cos[(m+2)x+\varphi_{m+2}] + \ldots +$$
$$+ \rho_{n-1}\cos[(n-1)x+\varphi_{n-1}]$$

By sum and dived by m-n:

$$V_m^{n-1}(f) = S_m + \sum_{k=m+1}^{n-1}\left(1 - \frac{k-m}{n-m}\right)\rho_k \cos(kx+\varphi_k)$$

**Remark (3):**

From Fourier series for function: $f(x) \approx \frac{a_0}{2} + \sum_{k=1}^{\infty}\alpha_k \cos kx + \beta_k \sin kx$

Therefore

$$f(x) - V_m^{n-1}(f) \approx \sum_{k=m+1}^{n-1}\left(\frac{k-m}{n-m}\right)\rho_k \cos(kx+\varphi_k) + \sum_{k=n}^{\infty}\rho_k \cos(kx+\varphi_k)$$

From Parseval equality

$$\|f(x)-V_m^{n-1}(f)\|^2 \approx \sum_{k=m+1}^{n-1}\left(\frac{k-m}{n-m}\right)^2 \rho_k^2 + \sum_{k=n}^{\infty}\rho_k^2 \qquad (4).$$

## 2. Main results

**Theorem (1):** for any function $f(x)\in L_2$, $(f(x)\neq const)$ and for any natural number $n$ will be obtained the inequality:

$$\|f(x)-\sigma_{n-1}(f)\|\leq \frac{1}{\sqrt{6}}\left\{\frac{n}{\pi}\int_0^{\frac{\pi}{n}}\omega_2^2(f,t)dt\right\}^{\frac{1}{2}} \qquad \forall n\geq n_0$$

**Proof:**

From Fourier series and Parseval equality will be obtained :

$$\|f(x-t)-2f(x)+f(x+t)\|^2 = 4\sum_{k=1}^{\infty}\rho_k^2(1-\cos kt)^2$$

On the other hand, with the define of second modulus of continuity

$$\omega_2(f,\delta)=\sup_{|t|\leq \delta}\|f(x-t)-2f(x)+f(x+t)\|$$

$$\omega_2^2(f,\delta)\geq \|f(x-\delta)-2f(x)+f(x+\delta)\|^2 \geq$$

$$\geq 4\sum_{k=1}^{\infty}\rho_k^2(1-\cos k\delta)^2$$

$$\geq 4\sum_{k=1}^{\infty}\rho_k^2(1-\cos k\delta)^2 \geq 4\sum_{k=1}^{n-1}\rho_k^2(1-\cos k\delta)^2 + 4\sum_{k=n}^{\infty}\rho_k^2(1-\cos k\delta)^2$$

$$\geq 4\sum_{k=1}^{n-1}\rho_k^2(1-\cos k\delta)^2 + 2\sum_{k=n}^{\infty}\rho_k^2(3-4\cos k\delta + \cos 2k\delta)$$

$$\frac{1}{6}\omega_2^2(f,\delta)\geq \frac{2}{3}\sum_{k=1}^{n-1}\rho_k^2(1-\cos k\delta)^2 + \sum_{k=n}^{\infty}\rho_k^2 - \frac{4}{3}\sum_{k=n}^{\infty}\rho_k^2\cos k\delta + \frac{1}{3}\sum_{k=n}^{\infty}\rho_k^2\cos 2k\delta$$

Changing $\delta$ to $t$ we can write

$$\sum_{k=n}^{\infty}\rho_k^2 \le \frac{1}{6}\omega_2^2(f,t) - \frac{2}{3}\sum_{k=1}^{n-1}\rho_k^2(1-\cos kt)^2 + \frac{4}{3}\sum_{k=n}^{\infty}\rho_k^2\cos kt - \frac{1}{3}\sum_{k=n}^{\infty}\rho_k^2\cos 2kt$$

Summing both sides by $\sum_{k=1}^{n-1}\left(\frac{k}{n}\right)^2\rho_k^2$ we obtain:

$$\sum_{k=n}^{\infty}\rho_k^2 + \sum_{k=1}^{n-1}\left(\frac{k}{n}\right)^2\rho_k^2 \le \frac{1}{6}\omega_2^2(f,t) + \sum_{k=1}^{n-1}\left(\frac{k}{n}\right)^2\rho_k^2 - \frac{2}{3}\sum_{k=1}^{n-1}\rho_k^2(1-\cos kt)^2 + \frac{4}{3}\sum_{k=n}^{\infty}\rho_k^2\cos kt - \frac{1}{3}\sum_{k=n}^{\infty}\rho_k^2\cos 2kt$$

from (3) we obtain

$$\|f(x)-\sigma_{n-1}(f)\|^2 \le \frac{1}{6}\omega_2^2(f,t) + \sum_{k=1}^{n-1}\rho_k^2\left[\left(\frac{k}{n}\right)^2 - \frac{2}{3}(1-\cos kt)^2\right] + \frac{4}{3}\sum_{k=n}^{\infty}\rho_k^2\cos kt - \frac{1}{3}\sum_{k=n}^{\infty}\rho_k^2\cos 2kt$$

in particular $\sum_{k=1}^{n-1}\rho_k^2\left(\frac{k}{n}\right)^2 \le \sum_{k=1}^{n-1}\rho_k^2\sin^2\left(\frac{kt}{2}\right)$ (see [1])

and hence

$$\|f(x)-\sigma_{n-1}(f)\| \le \frac{1}{6}\omega_2^2(f,t) + \sum_{k=1}^{n-1}\rho_k^2\left[\sin^2\left(\frac{kt}{2}\right) - \frac{2}{3}(1-\cos kt)^2\right] + \frac{4}{3}\sum_{k=n}^{\infty}\rho_k^2\cos kt - \frac{1}{3}\sum_{k=n}^{\infty}\rho_k^2\cos 2kt$$

$$\le \frac{1}{6}\omega_2^2(f,t) + \sum_{k=1}^{n-1}\rho_k^2\left[\frac{1}{2}(1-\cos kt) - \frac{2}{3}(1-\cos kt)^2\right] + \frac{4}{3}\sum_{k=n}^{\infty}\rho_k^2\cos kt - \frac{1}{3}\sum_{k=n}^{\infty}\rho_k^2\cos 2kt$$

$$\le \frac{1}{6}\omega_2^2(f,t) + \sum_{k=1}^{n-1}\rho_k^2\left[(1-\cos kt)\left(\frac{1}{2} - \frac{2}{3}(1-\cos kt)\right)\right] + \frac{4}{3}\sum_{k=n}^{\infty}\rho_k^2\cos kt - \frac{1}{3}\sum_{k=n}^{\infty}\rho_k^2\cos 2kt$$

$$\le \frac{1}{6}\omega_2^2(f,t) + \sum_{k=1}^{n-1}\rho_k^2\left[(1-\cos kt)\left(-\frac{1}{6} + \frac{2}{3}\cos kt\right)\right] + \frac{4}{3}\sum_{k=n}^{\infty}\rho_k^2\cos kt - \frac{1}{3}\sum_{k=n}^{\infty}\rho_k^2\cos 2kt$$

$$\le \frac{1}{6}\omega_2^2(f,t) + \frac{4}{3}\sum_{k=n}^{\infty}\rho_k^2\cos kt - \frac{1}{3}\sum_{k=n}^{\infty}\rho_k^2\cos 2kt$$

By integrating with respect to $t$ form $0$ to $\frac{\pi}{n}$ we obtain

$$\frac{\pi}{n}\|f(x)-\sigma_{n-1}(f)\|^2 \le \frac{1}{6}\int_0^{\frac{\pi}{n}}\omega_2^2(f,t)\,dt + \frac{4}{3}\sum_{k=n}^{\infty}\frac{\rho_k^2}{k}\sin\frac{k\pi}{n} - \frac{1}{6}\sum_{k=n}^{\infty}\frac{\rho_k^2}{k}\sin\frac{2k\pi}{n}$$

Let's prove $\sum_{k=n}^{\infty}\frac{\rho_k^2}{k}\sin\left(\frac{k\pi}{n}\right) \le 0$

We collect the terms of this series in the following form

$$\sum_{k=n}^{\infty} \frac{\rho_k^{\,2}}{k} \sin\left(\frac{k\pi}{n}\right) =$$

$$\sum_{l=1}^{\infty} \sum_{j=0}^{n-1} \left\{ \frac{\rho^{2}_{(2l-1)n+j}}{(2l-1)n+j} \sin[(2l-1)n+j]\frac{\pi}{n} + \frac{\rho^{2}_{2nl+j}}{2nl+j} \sin[2nl+j]\frac{\pi}{n} \right\}$$

On the other hand

$$\sin[(2l-1)n+j]\frac{\pi}{n} = \sin\left(2l\pi - \pi + \frac{j\pi}{n}\right) = \sin\left(\frac{j\pi}{n} - \pi\right) = -\sin\left(\frac{j\pi}{n}\right),$$

$$\sin[2nl+j]\frac{\pi}{n} = \sin\left(2l\pi + \frac{j\pi}{n}\right) = \sin\left(\frac{j\pi}{n}\right)$$

We obtain

$$\sum_{k=n}^{\infty} \frac{\rho^2}{k} \sin\left(\frac{k\pi}{n}\right) = \sum_{l=1}^{\infty} \sum_{j=0}^{n-1} \left[ \frac{\rho^{2}_{2nl+j}}{2nl+j} - \frac{\rho^{2}_{(2l-1)n+j}}{(2l-1)n+j} \right] \sin\left(\frac{j\pi}{n}\right)$$

We have

$$\forall l \geq 1 \Rightarrow 2nl+j > (2l-1)n+j$$

From remark (1)

$$\frac{\rho^2_{2nl+j}}{2nl+j} \leq \frac{\rho^2_{(2l-1)n+j}}{(2l-1)n+j} \qquad \forall l \geq 1 \quad , \quad \forall j \geq 0$$

And also

$$\sin\left(\frac{j\pi}{n}\right) \geq 0 \quad , \qquad 0 \leq j \leq n-1$$

Therefore

$$\sum_{l=1}^{\infty} \sum_{j=0}^{n-1} \left[ \frac{\rho^{2}_{2nl+j}}{2nl+j} - \frac{\rho^{2}_{(2l-1)n+j}}{(2l-1)n+j} \right] \sin\left(\frac{j\pi}{n}\right) \leq 0$$

$$\Rightarrow \sum_{k=n}^{\infty} \frac{\rho^2}{k} \sin\left(\frac{k\pi}{n}\right) \leq 0$$

Let's prove $\quad \sum_{k=n}^{\infty} \frac{\rho_k^{\,2}}{k} \sin\frac{2k\pi}{n} \geq 0$

To prove this inequality, we distinguish two cases

First case: Let $n$ be an odd number, then using the symbol $[\frac{n}{2}]$ To denote the correct division of a number $\frac{n}{2}$.

We collect the terms of this series in the following form

$$\sum_{k=n}^{\infty}\frac{\rho_k^2}{k}\sin\frac{2k\pi}{n} = \sum_{l=1}^{\infty}\{\sum_{j=0}^{[\frac{n}{2}]}\frac{\rho^2_{ln+j}}{ln+j}\sin 2(ln+j)\frac{\pi}{n} + \sum_{j=1}^{[\frac{n}{2}]}\frac{\rho^2_{(l+1)n-j}}{(l+1)n-j}\sin 2[(l+1)n-j]\frac{\pi}{n}\}$$

$$=\sum_{l=1}^{\infty}\{\sum_{j=0}^{[\frac{n}{2}]}\frac{\rho^2_{ln+j}}{ln+j}\sin\frac{2j\pi}{n} - \sum_{j=1}^{[\frac{n}{2}]}\frac{\rho^2_{(l+1)n-j}}{(l+1)n-j}\sin\frac{2j\pi}{n}\}$$

When $j=0$ we obtain $\sin\frac{2j\pi}{n}=0$ therefore

$$\sum_{k=n}^{\infty}\frac{\rho^2_k}{k}\sin\frac{2k\pi}{n} = \sum_{l=1}^{\infty}\sum_{j=1}^{[\frac{n}{2}]}(\frac{\rho^2_{ln+j}}{ln+j} - \frac{\rho^2_{(l+1)n-j}}{(l+1)n-j})\sin\frac{2j\pi}{n}$$

When $1\leq j \leq [\frac{n}{2}]$ we obtain $\sin\frac{2j\pi}{n} > 0$.

and $\quad\forall l \geq 1 \quad\quad ln+j < (l+1)n-j$

From remark (1)

$$\frac{\rho^2_{ln+j}}{ln+j} \geq \frac{\rho^2_{(l+1)n-j}}{(l+1)n-j}$$

Therefore

$$\sum_{l=1}^{\infty}\sum_{j=1}^{[\frac{n}{2}]}(\frac{\rho^2_{ln+j}}{ln+j} - \frac{\rho^2_{(l+1)n-j}}{(l+1)n-j})\sin\frac{2j\pi}{n} \geq 0$$

$$\Rightarrow \sum_{k=n}^{\infty}\frac{\rho^2_k}{k}\sin\frac{2k\pi}{n} \geq 0$$

Second case: Let $n$ be an even number

We collect the terms of this series in the following form

$$\sum_{k=n}^{\infty}\frac{\rho_k}{k}\sin\frac{2k\pi}{n}=\sum_{l=1}^{\infty}\{\sum_{j=0}^{\frac{n}{2}-1}\frac{\rho^2_{ln+j}}{ln+j}\sin 2(ln+j)\frac{\pi}{n}+\sum_{j=1}^{\frac{n}{2}}\frac{\rho^2_{(l+1)n-j}}{(l+1)n-j}\sin 2[(l+1)-j]\frac{\pi}{n}\}=$$

$$=\sum_{l=1}^{\infty}\{\sum_{j=0}^{\frac{n}{2}-1}\frac{\rho^2_{ln+j}}{ln+j}\sin\frac{2j\pi}{n}-\sum_{j=1}^{\frac{n}{2}}\frac{\rho^2_{(l+1)n-j}}{(l+1)n-j}\sin\frac{2j\pi}{n}\}$$

When $j=0$ we obtain $\sin\frac{2j\pi}{n}=0$

When $j=\frac{n}{2}$ we obtain $\sin\frac{2j\pi}{n}=0$ therefore

$$\sum_{k=n}^{\infty}\frac{\rho^2_k}{k}\sin\frac{2k\pi}{n}=\sum_{l=1}^{\infty}\sum_{j=1}^{\frac{n}{2}-1}(\frac{\rho^2_{ln+j}}{ln+j}-\frac{\rho^2_{(l+1)n-j}}{(l+1)n-j})\sin\frac{2j\pi}{n}\geq 0$$

Therefore

$$\frac{\pi}{n}\|f(x)-\sigma_{n-1}(f)\|^2\leq\frac{1}{6}\int_0^{\frac{\pi}{n}}\omega_2^{\,2}(f,t)\,dt$$

$$\Rightarrow\|f(x)-\sigma_{n-1}(f)\|\leq\frac{1}{\sqrt{6}}\left\{\frac{n}{\pi}\int_0^{\frac{\pi}{n}}\omega_2^{\,2}(f,t)\,dt\right\}^{\frac{1}{2}}$$

**Theorem (2):** for any function $f(x)\in L_2$, $(f(x)\neq const)$ and for any natural number $m<n$ ; $n,m$ Will be obtained the inequality:

$$\|f(x)-V_m^{n-1}(f)\|\leq\frac{1}{\sqrt{6}}\left\{\frac{n}{\pi}\int_0^{\frac{\pi}{n}}\omega_2^{\,2}(f,t)\,dt\right\}^{\frac{1}{2}} \qquad \forall n\geq n_0$$

**Proof:**

To prove the inequality it's enough to prove

$$\|f(x)-V_m^{n-1}(f)\|\leq\|f(x)-\sigma_{n-1}(f)\|$$

We have

$$\frac{k-m}{n-m}\leq\frac{k}{n}$$

Therefore

$$\sum_{k=m+1}^{n-1}\left(\frac{k-m}{n-m}\right)^2\rho^2 \leq \sum_{k=m+1}^{n-1}\left(\frac{k}{n}\right)^2\rho^2$$

Therefore

$$\sum_{k=m+1}^{n-1}\left(\frac{k-m}{n-m}\right)^2\rho^2 + \sum_{k=n}^{\infty}\rho_k^2 \leq \sum_{k=m+1}^{n-1}\left(\frac{k}{n}\right)^2\rho^2 + \sum_{k=n}^{\infty}\rho_k^2$$

$$\left\|f(x)-V_m^{n-1}(f)\right\|^2 = \sum_{k=m+1}^{n-1}\left(\frac{k-m}{n-m}\right)^2\rho^2 + \sum_{k=n}^{\infty}\rho_k^2 \leq \sum_{k=m+1}^{n-1}\left(\frac{k}{n}\right)^2\rho^2 + \sum_{k=n}^{\infty}\rho_k^2 = \left\|f(x)-\sigma_{n-1}(f)\right\|^2$$